\documentclass{icmart}

\contact[vergne@math.polytechnique.fr]{Mich{\`e}le Vergne;
Centre de Math{\'e}matiques Laurent Schwartz, 91128 Palaiseau, France} 

\usepackage{amssymb,amsmath}

\usepackage{graphicx}

\newtheorem{theorem}{Theorem}
\newtheorem{proposition}[theorem]{Proposition}
\newtheorem{definition}[theorem]{Definition}
\newtheorem{corollary}[theorem]{Corollary}
\newtheorem{example}[theorem]{Example}
\newtheorem{remark}[theorem]{Remark}
\renewcommand{\t}{{\mathfrak t}}
\newcommand{\g}{{\mathfrak g}}
\newcommand{\h}{{\mathfrak h}}
\renewcommand{\k}{{\mathfrak k}}
\renewcommand{\r}{{\mathfrak \r}}
\newcommand{\C}{{\mathbb C}}
\newcommand{\R}{{\mathbb R}}
\newcommand{\Z}{{\mathbb Z}}

\newcommand{\CA}{{\cal A}}
\newcommand{\CB}{{\cal B}}

\newcommand{\CE}{{\cal E}}
\newcommand{\CF}{{\cal F}}
\newcommand{\CH}{{\cal H}}
\newcommand{\CL}{{\cal L}}

\newcommand{\vol}{\operatorname{vol}}
\newcommand{\Tr}{\operatorname{Tr}}
\newcommand{\form}{{\alpha}}
\newcommand{\res}{\operatorname{res}}
\newcommand{\Dirac}{{\delta_0}}

\newcommand{\bx}{\mathbf{x}}

\newcommand{\gc}{\mathfrak{c}}
\newcommand{\beps}{\boldsymbol{\epsilon}}

\title{Applications of Equivariant Cohomology}

\author{Mich\`ele Vergne}

\date{}
\begin{document}

\begin{abstract}
We will discuss the equivariant cohomology of a manifold endowed
with the action of a Lie group. Localization formulae for
equivariant integrals are explained by a vanishing theorem for
equivariant cohomology with generalized coefficients. We then give
applications to integration of characteristic classes on
symplectic quotients and to indices of transversally elliptic
operators. In particular, we state a conjecture for the index of a
transversally elliptic operator linked to  a Hamiltonian action.
In the last part, we describe algorithms for numerical
computations of values of multivariate spline functions and  of
vector-partition functions of classical root systems.
\end{abstract}

\begin{classification}
Primary: $K$-theory; Secondary: Convex and Discrete geometry.
\end{classification}

\begin{keywords}
equivariant cohomology, Hamiltonian action, symplectic reduction,
localization formula, polytope, index, transversally elliptic
operator, spline, Euler-Maclaurin formula.
\end{keywords}

\maketitle


\section{Introduction}\label{intro}

The aim of this article is to show  how theorems of localization
in equivariant cohomology   not only provide beautiful
mathematical formulae, but also stimulated progress in algorithmic
computations. I will focus on my favorite themes: quantization of
symplectic manifolds and algorithms for polytopes, and neglect
many other applications. Many mathematicians have shared their
ideas with me, notably  Welleda Baldoni, Nicole Berline, Michel
Brion, Michel Duflo, Shrawan Kumar, Paul-Emile Paradan and Andras
Szenes.  I will therefore  often employ a collective ``we",
instead of anxiously weighing my own contribution.

I will describe here the theory of equivariant cohomology with
generalized coefficients of a manifold $M$ on which a Lie group
$K$ acts. The integral  of such a cohomology class is a
generalized function  $I(\phi)$ on $\k$, with $\phi$ in $\k$, the
Lie algebra of $K$. We wish to solve  two problems. The first  is
to give a ``localization formula" for $I(\phi)$ as a ``short"
expression. The second  is: given such a short formula for
$I(\phi)$, compute the value $\hat I(\xi)$ of the Fourier
transform of $I$ at a point $\xi\in \k^*$ in terms of  the initial
geometric data. Let me give the motivation for  such questions.

By integrating  de Rham cohomology classes on a manifold, one
obtains certain numerical quantities. For example, the symplectic
volume $\vol_M$ of a compact symplectic manifold $M$ is the
integral of the Liouville form, and the Atiyah-Singer
cohomological formula for the index of an elliptic operator $D$ on
$M$  is  an integral of a cohomology class with compact support on
$T^*M$. In the interplay between toric varieties and polytopes,
these numerical quantities correspond respectively to the volume
of a polytope and to the number of integral points inside a
rational polytope. Moreover, the volume is the classical limit of
the discrete version, the number of points in dilated polytopes.

 When the
manifold is provided with the action of a compact Lie group $K$,
similar objects are described by integrals of equivariant
cohomology classes. The equivariant volume $\vol_M(\phi)$ of a
compact Hamiltonian manifold $M$ is a $C^{\infty}$ function of
$\phi\in \k$, obtained by integrating  a particular  equivariant
cohomology class on $M$. More generally, if  $M$ is a non-compact
Hamiltonian manifold with proper moment map, and additional
convergence conditions, its equivariant volume $\vol_M(\phi)$ is
  a generalized function on $\k$. As shown by Duistermaat-Heckman, the
value at  $\xi\in \k^*$ of the inverse Fourier transform  of
$\vol_M(\phi)$ is the symplectic volume of the Marsden-Weinstein
reduction of $M$ at $\xi$. If a $K$-invariant operator $D$ is
elliptic in the directions transverse to the orbits of $K$, its
index ${\rm Index}(D)$  is a generalized function on $K$, that is,
a series of characters of $K$. It can be described  in terms of
integrals of equivariant cohomology classes on $T^*M$. The
discrete inversion problem is to determine  each Fourier
coefficient of ${\rm Index}(D)$. When $D$ is an operator linked to
the symplectic structure, we  think of ${\rm Index}(D)$ as the
quantum version of the equivariant volume. The Guillemin-Sternberg
conjecture, now established by Meinrenken-Sjamaar for any compact
Hamiltonian manifold, is an example where such an inversion
problem  has a beautiful answer in geometric terms.

In the case of a manifold with a  circular symmetry, we proved a
localization formula for integrals of equivariant cohomology
classes as a sum of local contributions from the fixed points.
This formula
 is  similar to the Atiyah-Bott Lefschetz fixed point formula
for the equivariant index of an elliptic operator on $M$. A
drawback of such formulae is that each individual term has poles,
and the Fourier transform of an individual term is  meaningless.
We will describe here a more general principle of localization for
integrals of equivariant cohomology classes. Let $\kappa$ be a
$K$-invariant vector field tangent to the orbits of $K$. Witten
showed that equivariant integrals on $M$ can be computed in terms
of local data near the set $C$ of zeroes of $\kappa$. Furthermore,
for each connected component $C_F$ of $C$, the local contribution
of $C_F$ is a generalized function on $\k$. Witten's localization
theorem can be best understood through Paradan's identity: $1=0$
on $M-C$, in equivariant cohomology with generalized coefficients.
Basic definitions and Paradan's identity are explained in Section
\ref{equi}.

The identity $1=0$ on $M-C$  has many independent  applications
that we describe in Section \ref{appli}. When $M$ is a Hamiltonian
manifold with moment map $\mu$, the set of zeroes of the Kirwan
vector field is the set of critical points of the function
$\|\mu\|^2$. According to  Witten's theorem, integrals on reduced
spaces of $M$ can be related to equivariant integrals on $M$.
Using a similar localization argument for transversally elliptic
operators, Paradan was able to extend the proof of the
Guillemin-Sternberg conjecture to some non-compact Hamiltonian
spaces linked to representation theory of real semi-simple Lie
groups via Kirillov's orbit method. We will state  a
generalization of the Guillemin-Sternberg conjecture for a
transversally elliptic operator canonically attached to a
Hamiltonian action in Section \ref{appli}.

From the localization formulae, one is led to study  generalized
functions which are regular outside a union of hyperplanes. This
will be the topic of Section \ref{arrangement}. In particular, we
will relate the cohomology ring of toric manifolds to cycles in
the complement of an arrangement of hyperplanes.

 As there are some relations between Hamiltonian
geometry and convex polytopes, these localization theorems have an
analogue for polytopes. Such an analogue is the local
Euler-Maclaurin formula for polytopes, which was conjectured by
Barvinok-Pommersheim. We will  indicate in Section \ref{poly} how
some theoretical results on intersection rings can be turned
 into effective tools for
numerical computations. We implemented algorithms  for various
problems such as computing the value of the convolution of a large
number of Heaviside distributions, the number of integral points
in network polytopes and   Kostant partitions functions, with
applications to the tensor multiplicities formulae. This last
section  can be read independently. Indeed, these applications to
polytopes have elementary proofs, but it was through interaction
with Hamiltonian geometry that some of these tools were
discovered.

For lack of space, I was only able to include  central references
to the topics discussed in this text. For more bibliographical
comments, references and motivations, one might consult
\cite{ber-get-ver}, \cite{duf-kum-ver}, \cite{gui-ste99},
\cite{ver01} and my home page (notably, the text called
``Ex{\'e}g{\`e}se") at math.polytechnique.fr/cmat/vergne/ . The texts
\cite{vergnewomen} and \cite{ver99} are introductory and hopefully
reader-friendly.

\section{Simple Examples}\label{baby}

In this section, I will give simple examples of sums  which can be
represented by  short formulae, and  a simple example of the
inverse problem we have in mind. I will also give  a  sketch of
the proof of  the  stationary phase formula  as similar stationary
phase arguments will be our fundamental tools.

\subsection{Geometric series}\label{geo}
Some formulae in mathematics condense a large amount of
information in short expressions. The most striking formula
perhaps is the one that sums a very long geometric series:

$$\sum_{i=0}^{10000}q^i=\frac{1}{1-q}+\frac{q^{10000}}{1-q^{-1}}.$$

For a  straightforward calculation of the left hand side for a
given value $q$, one needs  to know  the value of the function
$q^i$ at all the $10001$ integral points of the interval
$[0,10000]$, while for the right hand side one needs only the
value of this function at the end points $0,10000$. Note that each
term of the right hand side has a  pole  at $q=1$.

\vspace{1cm}

\includegraphics{mesdebuts.05}

\vspace{1cm}

The short formula (here $A,B,i$ are integers)
\begin{equation}\label{shortgeo}
\sum_{i=A}^B
q^i=\frac{q^A}{1-q}+\frac{q^B}{1-q^{-1}}=-\frac{q^{A-1}}{1-q^{-1}}-\frac{q^{B+1}}{1-q}
\end{equation}
is related to the  following equalities   of characteristic
functions:
\begin{eqnarray*}
\chi([A,B])&=&\chi([A,\infty[)+\chi(]-\infty,B])-\chi(\R)\\
 &=&\chi(\R)-\chi(]-\infty,A[)-\chi(]B,\infty[). \end{eqnarray*}

We draw the picture of the last equality.

\begin{figure}[!h]
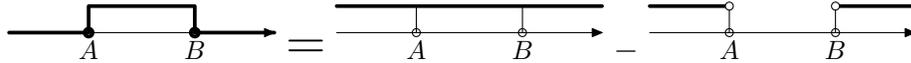

\includegraphics{decomp.01}
{\huge$\mathbf{ =}$}
\includegraphics{decomp.02}
  {\large$\mathbf {-}$}
\includegraphics{decomp.03}
\caption{Decomposition of an interval}\label{decomint}
\end{figure}

 Then to sum  $q^i$ from $A$ to $B$, we first sum $q^i$
from $-\infty$ to $\infty$ and subtract the two sums over the
integers strictly less  than $A$ and over the integers strictly
greater than $B$. Thus, if
$$S_0:=\sum_{i=-\infty}^{\infty}
q^i,\hspace{1cm}S_A:=\sum_{-\infty}^{A-1} q^i, \hspace{1cm}
S_B:=\sum_{B+1}^{\infty}q^i,$$ we obtain formally, or, setting
$q=e^{i\phi}$, in the sense of generalized functions on the unit
circle,
\begin{equation}\label{geotrans}
S=S_0-S_A-S_B.
\end{equation}

For a  value $q\neq 1$, the first sum $S_0$ is $0$ as follows from
$(1-q)S_0=0$, while $S_A$, $S_B$ are just geometric progressions
and we come back to the short formula (\ref{shortgeo}).

The reader may recognize in  Formula (\ref{shortgeo}) a very
simple instance of the Atiyah-Bott Lefschetz fixed point formula
on the Riemann sphere. Formula (\ref{geotrans}) illustrates
Paradan's localization of elliptic operators, which we describe in
Section \ref{transver}. Indeed, Formula (\ref{geotrans}) is an
example of the decomposition of the equivariant index of an
elliptic operator on the Riemann sphere in a sum of indices of
three transversally elliptic operators (see Example
\ref{ellipticonsphere}).

\subsection{Inverse problem}\label{inverse}
The inverse problem may be described as follows: given a short
expression for a sum, compute an individual term of the sum.

 Here is an example. Consider the
following product of  geometric series $G:=(\sum_{i=0}^{\infty}
q_1^i)^3(\sum_{j=0}^{\infty} q_2^j)^3(\sum_{k=0}^{\infty}
q_1^kq_2^k)^3$ given by the short expression:
$$S(q_1,q_2):=\frac{1}{(1-q_1)^3}\frac{1}{(1-q_2)^3}\frac{1}{(1-q_1q_2)^3}.$$
Let us  compute the coefficient $c(a,b)$ of $q_1^{a}q_2^{b}$ in
$G$. If $a\geq b$, an iterated application of the residue theorem
in one variable leads to
$$c(a,b)=\res_{x_2=0}\left(\res_{x_1=0}
\frac{e^{ax_1}e^{bx_2}\, dx_1\,
dx_2}{(1-e^{-x_1})^3(1-e^{-x_2})^3(1-e^{-(x_1+x_2)})^3}\right).$$
If we set
$$g(a,b)=\frac{(b+1)(b+2)(b+3)(b+4)(b+5)(7 a^2-7 ab +2b^2+21
a-9b+14)}{14\cdot 5! },$$  we obtain the following equalities.
\begin{eqnarray}\label{factor}
\text{If}\, a\geq b, \,\,\text {then} \, c(a,b)&=&g(a,b).\\
\text{If}\, a\leq b, \,\,\text {then} \, c(a,b)&=&g(b,a).
\end{eqnarray}

We will discuss in Section \ref{Heaviside} a residue theorem
(Theorem \ref{cycle}) in several variables, which gives an
algorithmic solution to this type of inversion problem.

The Guillemin-Sternberg conjecture (see Section \ref{QR}) gives a
geometric interpretation of the Fourier coefficients of series for
similar inversion problems.

\subsection{Stationary phase}\label{station}

Let $M$ be a compact manifold of dimension $n$, $f$ a smooth
function on $M$ and $dm$ a smooth density. Consider the function
$$F(t):=\int_{ M}e^{itf(m)}dm.$$
The dominant contribution to the value of this integral as $t$
tends to infinity arises from the neighborhood of the set $C$ of
critical points of $f$. We indicate a proof of this fact, as
similar arguments will be employed later on. Consider the image of
$M$ by the map $x=f(m)$ and the push-forward of the density $dm$.
Then $F(t)=\int_{\R}e^{it x}f_{*}(dm).$ Choose a smooth function
$\chi$ on $M$, equal to $1$ in a neighborhood of the set $C$ and
supported near $C$. Then $F(t)=F_C(t)+R(t),$ where
$$F_C(t):=\int_{\R}e^{it x} f_*(\chi dm),\hspace{1cm}
R(t)=\int_{\R}e^{it x} f_*((1-\chi) dm).$$ $R(t)$ is the Fourier
transform of a smooth compactly supported function, and thus
decreases rapidly at infinity. It is not hard to show that, if $f$
has a finite number of non-degenerate critical points, then
$$F(t)\sim F_C(t)\sim \sum_{p\in C}e^{itf(p)}\sum_{k\geq 0}^{\infty} a_{p,k}
t^{-\frac{n}{2}+k},$$ where the constants $a_{p,k}$ can be
computed in terms of $f$ and $dm$ near $p\in C$. We can say that
asymptotically, the integral ``localizes" at a finite number of
points $p$.

\begin{figure}[!h]
\begin{center}\includegraphics{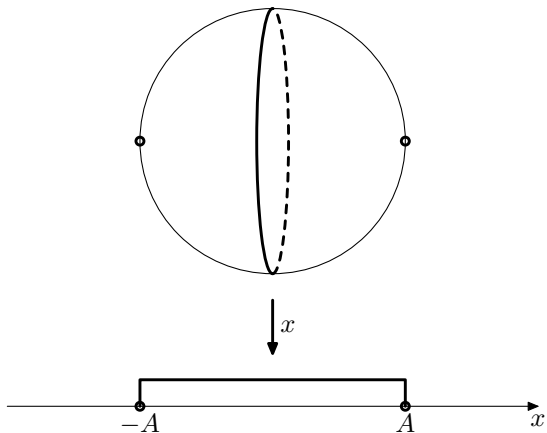}
\end{center}\caption{Projecting the sphere $x^2+y^2+z^2=A^2$}\label{sphere}
\end{figure}

\begin{example}\label{stationaryonsphere}
  Let $M$ be the sphere $\{x^2+y^2+z^2=A^2\}$ of radius $A$ endowed
  with the Liouville volume form $dm:=\frac{dy\wedge dz}{2\pi x}$. Let our function $f$ be the projection onto the $x$-axis: $f=x$. We
   immediately see that $f_*(dm)$ is the characteristic function of
  the interval $[-A,A]$. Thus we obtain  the formula:
 \[
F(t)=\int_{-A}^A
e^{itx}dx=\frac{e^{-iAt}}{-it}+\frac{e^{iAt}}{it}.\]
\end{example}

Observe that here $F(t)$ is not just asymptotically, but exactly
equal to the local expression. The reason is that in this example
the function $f$ is the Hamiltonian of an action of the circle
group $S^1:=\{e^{i\phi}\}$ on a compact symplectic manifold, and
$dm$ is the Liouville measure. In such a case, the
Duistermaat-Heckman exact stationary phase formula \cite{dui-hec}
implies that $f_*(dm)$ is locally polynomial on $f(M)$ and that
\begin{equation}\label{stationaryphase}
F(t)=\sum_{p\in C}e^{itf(p)} a_{p,0} t^{-\frac{n}{2}}.
\end{equation}

We will interpret the Duistermaat-Heckman formula as an example of
the abelian localization formula (Theorem \ref{BVAB}) of integrals
of equivariant forms in Section \ref{wit}.

\section{Equivariant differential forms}\label{equi}

Our motivation to study equivariant differential forms came from
representation theory.

 Let $M$ be a manifold with an action of the
circle group $S^1$. The Atiyah-Bott fixed point formula
\cite{ati-bot67} describes the equivariant index of an elliptic
operator on $M$ in terms of local data near the fixed points of
the action. One of the applications of the formula was a geometric
interpretation of the Weyl formula for the characters of
irreducible representations of compact Lie groups.

The character formula has continuous analogues: the formulae for
the Fourier transforms of coadjoint orbits, which are linked to
representation theory via Kirillov's orbit method. For compact
groups this is the Harish-Chandra formula; for non-compact
semi-simple groups, Rossmann gave a fixed point formula in the
case of discrete series characters.

In joint work with Nicole Berline, I found a geometric
interpretation of Rossmann's formula using equivariant forms
\cite{ber-ver82-1}. The cohomological tool behind our computation
was a deformation of the de Rham complex with the use of vector
fields. A similar formalism was described by Witten \cite{wit82}
with different motivation. There were earlier results which
condensed certain integrals on $M$ in short formulae localized
near ``fixed points", such as Bott's  residue formulae
\cite{Bot67}, its generalization by Baum-Cheeger  \cite{bau-che}
and the Duistermaat-Heckman exact stationary phase formula
\cite{dui-hec}. As explained by Atiyah and Bott \cite{ati-bot84},
our result was related to localization in topological equivariant
cohomology. However, this revival of ``de Rham" theory of
equivariant cohomology in  terms of differential forms turned out
to be very fruitful, especially in applications to non-compact
spaces and stationary phase type arguments.

\subsection{Equivariant de Rham complex}\label{forms}

{\bf Notation}: I  keep the notation $N$ for not necessarily
compact manifolds, and $M$ for compact manifolds. Similarly a
compact group will be denoted by the letter $K$, while $G$ will be
an arbitrary real Lie group. The letters $T,H$ will be reserved
for  tori, which are compact connected abelian Lie groups, and
therefore are just products of circle groups $\{e^{i \theta_a}\}$.
In this case, I take as basis of the Lie algebra $\t$, elements
$J_a$ such that $\exp(\theta_a J_a):=e^{i\theta_a}$ $(\theta_a\in
\R)$. The gothic german letters $\g, \k, \t, \h$  denote the
corresponding Lie algebras, $\g^*, \k^*, \t^*, \h^*$ the dual
vector spaces, $J^a$ the dual basis to a basis $J_a$. If $s\in G$,
I denote by $N_s$ the set of fixed points of the action of $s$ on
the $G$-manifold $N$. The letter $\phi$ denotes  an element of
$\g$. If $\g=\R J$ is the Lie algebra of $S^1$, I identify $\g$
and $\R$. I denote by $S(\g^*)$ the algebra of polynomial
functions on $\g$, by $C^{\infty}(\g)$ the space of $C^{\infty}$
functions on $\g$ and by $C^{-\infty}(\g)$ the space of
generalized functions on $\g$. An element $v\in C^{-\infty}(\g)$
is denoted by $v(\phi)$ although the value at $\phi\in \g$ of $v$
may not be defined. By definition, it is always defined in the
distributional sense: if $F(\phi)$ is a $C^{\infty}$ function on
$\g$ with compact support (a test function), then $\langle v,F
d\phi\rangle$, denoted by $\int_\g v(\phi)F(\phi)d\phi$, is well
defined.

\bigskip

Let us first define the  equivariant cohomology algebra with
$C^{\infty}$ coefficients of a $G$-manifold $N$.

Let $G$ be a Lie group acting on a manifold $N$. For $\phi\in \g$,
we denote by $ V\phi$ the vector field on $N$ generated by the
infinitesimal action of $-\phi$ :  for $x\in N$,
$V_x\phi:=\frac{d}{d\epsilon} \exp(-\epsilon \phi)\cdot
x|_{\epsilon=0}.$ If $N$ is provided with an action of $S^1$, we
simply denote by $J$ the vector field $VJ$. Let $\CA(N)$ be the
algebra of differential forms on $N$ with complex coefficients,
and denote by $d$ the exterior derivative. If $V$ is a vector
field, let $\iota(V)$ be the contraction by $V$. If
$\nu:=\sum_{i=0}^{\dim N}\nu_{[i]}$ is a differential form on an
oriented manifold $N$, then the integral of $\nu$ over $N$ is by
definition the integral of the top degree term of $\nu$: $\int_N
\nu:=\int_N \nu_{[\dim N]}$, provided that this last integral is
convergent.

A smooth map $\form:\g\to \CA(N)$ is called an {\em equivariant
form}, if $\form$ commutes with the action of $G$ on both sides.
The equivariant de Rham operator $D$ (\cite{ber-ver82-1},
\cite{wit82}) may be viewed as a deformation of the de Rham
operator $d$ with the help of the vector field $V\phi$. It is
defined on equivariant forms by the formula
$$(D(\form))(\phi):=d(\form(\phi))-\iota(V\phi)\form(\phi).$$
Then $D^2=0$. An equivariant form $\form$ is equivariantly closed
if $D\form=0$. The {\em cohomology space}, denoted by
$\CH^{\infty}(\g,N)$, is, as usual, the kernel of $D$ modulo its
image. This is an  algebra, $\Z/2\Z$-graded in even and odd
classes. If $G:=\{1\}$, this is the usual cohomology algebra
$\CH(N)$.

The integral of an equivariant differential form  may be defined
as a generalized function. Indeed, let $F(\phi)$ be a test
function on $\g$; then $\int_\g \form(\phi)F(\phi)d\phi$ is a
differential form on $N$. If this differential form is integrable
on $N$ for all test functions $F$, then $\int_N\alpha$ is defined
by
$$\langle \int_N \form,Fd\phi\rangle :=\int_N \int_\g \form(\phi) F(\phi)d\phi.$$

Of course if $N$ is compact oriented, $\int_N\form(\phi)$ is a
$C^{\infty}$ function.

\subsection{Hamiltonian spaces}\label{Hamil}
Examples of equivariantly closed forms  arise in Hamiltonian
geometry.

Let $N$ be a symplectic manifold with symplectic form $\Omega$. We
say that the action of $G$ on $N$ is Hamiltonian with moment map
${\bf \mu}:N\to \g^*$ if, for every $\phi\in \g$,\, $d(\langle
\phi,{\bf \mu}\rangle )=\iota(V\phi)\cdot \Omega$. Thus the zeroes
of the vector field $V\phi$ (that is, the {\bf fixed points} of
the one parameter group $\exp(t\phi)$) are the critical points of
$\langle\phi,\mu\rangle $.

The equivariant symplectic form
 $\Omega(\phi):=\langle \phi,{\bf \mu}\rangle +\Omega$\, is a closed
equivariant form. Indeed,  $$(d-\iota(V\phi))(\langle \phi,{\bf
\mu}\rangle +\Omega)=d(\langle \phi,{\bf \mu}\rangle
)-\iota(V\phi)\cdot \Omega+d(\Omega)$$ and this is equal to  $0$
as both equations
$$
d\Omega=0, \hspace{1cm} d(\langle \phi,{\bf \mu}\rangle
)=\iota(V\phi)\cdot \Omega$$ hold.

The two basic examples of  Hamiltonian spaces with an Hamiltonian
action of $S^1$ are:

{$(1)$} $\R^2$ if the action of $S^1$ has a fixed point.

{$(2)$} the cotangent bundle $T^*S^1$ if the action of $S^1$ is
free.

\bigskip

{$(1)$} Let $N:=\R^2$ with coordinates $[x,y]$. The circle group
$S^1$ acts by rotations with  isolated fixed point $[0,0]$. The
symplectic form is $\Omega:=dx\wedge dy$. The function
$\frac{x^2+y^2}{2}$ is the Hamiltonian function for the vector
field $J:= y\partial_x -x\partial_y$. Thus the equivariant
symplectic form is
$$\Omega(\phi)=\phi\left(\frac{x^2+y^2}{2}\right)+ dx\wedge dy.$$

{$(2)$} Let $N:=T^*S^1=S^1\times \R$. The circle group $S^1$ acts
freely by rotations on $S^1$. If $[e^{i\theta},t]$ is a point of
$T^*S^1$ with $t\in \R$, the symplectic form is $\Omega:=dt\wedge
d\theta$. The function $t$ is the Hamiltonian function for the
vector field $J:=-
\partial_\theta$. Thus the equivariant symplectic form is
$$\Omega(\phi)= \phi \,t+ dt\wedge d\theta.$$\\

A particularly important closed equivariant form is
$e^{i\Omega(\phi)}$. If $\dim N:=2\ell$, then
$$e^{i\Omega(\phi)}=e^{i\langle \phi,{\bf \mu}\rangle }\left(1+i\Omega+\frac{(i\Omega)^2}{2!}+\cdots
+\frac{(i\Omega)^{\ell}}{\ell !}\right).$$

\subsection{Equivariant volumes}\label{equivol}

Let $M$ be a compact  $K$-Hamiltonian  manifold of dimension
$2\ell$. By definition, the equivariant symplectic volume of $M$
is the function of $\phi\in \k$ given by
$$\vol_M(\phi):=\frac{1}{(2i\pi)^{\ell}}\int_M e^{i\Omega(\phi)}=
\int_M e^{i\langle \phi,\mu(m)\rangle } \frac{\Omega^{\ell}}{\ell
!(2\pi)^\ell}.$$  Note that $\vol_M(0)$ is the symplectic volume
of $M$. The last integral, according to the
Duistermaat-Heckman-formula \cite{dui-hec}, localizes as a sum of
integrals on the connected components of the set of zeroes of
$V\phi$. If this set of zeroes is finite,
\begin{equation}\label{DH}
\vol_M(\phi)=\sum_{p\in \rm zeroes\, of\, V\phi}
\frac{e^{i\langle\phi,\mu(p)\rangle
}}{i^{\ell}\sqrt{\det_{T_pM}L_p(\phi)}}, \end{equation} where
$L_p(\phi)$ is the endomorphism of $T_pM$ determined by the
infinitesimal action of $\phi$ at $p$. As $K$ is  compact, there
is a well-defined polynomial square root of the function
$\phi\mapsto \det_{T_pM}L_p(\phi)$, the sign being determined by
the  orientation.

\begin{example}\label{stationaryonsphere2}
  Consider, as in Example \ref{stationaryonsphere}, Section \ref{station}, the sphere $M$
  with
$S^1$-action given by rotation around the $x$-axis and
$\Omega:=\frac{dy\wedge dz}{x}$. Then $f:=x$ is the Hamiltonian
function  of the vector field  $J:= (y\partial_z-z\partial_y)$.
The equivariant volume is the $C^{\infty}$ function
$$\vol_M(\phi)=\int_M e^{i\phi x} dm=\frac{e^{-iA\phi}}{-i\phi}+\frac{e^{iA
\phi}}{i\phi}.$$
\end{example}

\bigskip

Let us point out some examples of non-compact manifolds $N$ where
the equivariant symplectic volume exists in the sense of
generalized functions. We will use the following generalized
functions:
$$
Y^+(\phi):=\int_{0}^{\infty} e^{i\phi t}dt,\hspace{1cm}
Y^-(\phi):=\int_{-\infty}^{0} e^{i\phi t}dt, \hspace{1cm}
\Dirac(\phi):=\int_{-\infty}^{\infty} e^{i\phi t}dt.$$

Note that the generalized function $Y^+(\phi)$  is the boundary
value of the holomorphic function $\frac{1}{-i\phi}$ defined on
the upper-half plane, so that it satisfies the relation
$(-i\phi)Y^+(\phi)=1$. The generalized function $\Dirac(\phi)$
satisfies the relation $\phi\, \Dirac(\phi)=0$.

Return to our two basic examples $\R^2$ and $T^*S^1$ with action
of $S^1$.

{$(1)$} $N:=\R^2$. We have
\begin{equation}\label{2SD}
 \vol_N(\phi)=\frac{1}{2\pi}\int_{\R^2}e^{i\phi
\frac{(x^2+y^2)}{2}}dx dy=\int_{0}^{\infty} e^{i\phi
r}dr=Y^{+}(\phi).
\end{equation}
 When $\phi\neq
0$, we have $\vol_N(\phi)=\frac{1}{-i\phi}$. This coincides with
what would be the Duistermaat-Heckman formula in the non-compact
case: there is just one fixed point $[0,0]$ for the action.

$(2)$ $N:=T^*S^1$. We have
$$\vol_N(\phi)=\frac{1}{2\pi}\int_{\R\times S^1}e^{i  \phi t}dt d\theta=
\int_{\R}e^{i \phi t}dt= \Dirac( \phi).$$ Thus $\vol_N(\phi)$ is
always $0$ when $\phi\neq 0$. This is consistent with the fixed
point philosophy: the action of $S^1$ on $T^*S^1$ is free, thus
the set of zeroes of $V\phi$ is empty when $\phi\neq
0$.\\

The next example illustrates our original motivation to introduce
the equivariant differential complex.

{\bf Coadjoint orbits.} Let $G$ be a real Lie group. Recall
\cite{kos} that when $N:=G\lambda$ is the orbit of an element
$\lambda\in \g^*$ by the coadjoint representation, then $N$ has a
$G$-Hamiltonian structure, such that the moment map is the
inclusion $N\to \g^*$. The equivariant volume $\vol_N(\phi)$ is
defined as a generalized function on $\g$, if the orbit $G\lambda$
is tempered. This is just the Fourier transform of the
$G$-invariant measure supported on $G\lambda\subset \g^*$.

 When $N$ is a coadjoint orbit of a
compact Lie group $K$, Harish-Chandra  gave a fixed point formula
for $\vol_N(\phi)$. Now this is seen as a special case of the
Duistermaat-Heckman formula (\ref{DH}). Rossmann \cite{ros1} and
Libine \cite{lib04} extended the Harish-Chandra formula
 to the case of closed coadjoint orbits of reductive non-compact Lie groups, involving delicate constants at fixed points
at ``infinity" defined combinatorially by Harish-Chandra and Hirai
 and  topologically by Kashiwara.

Here is an example. Consider the group $SL(2,\R)$ with Lie algebra
 $\g$ with basis

$$J_1:=\left(
\begin{array}{cc}
  1 & 0 \\
  0 & -1 \\
\end{array}
\right),\hspace{1cm} J_2:=\left(
\begin{array}{cc}
  0 & 1 \\
  1 & 0 \\
\end{array}
\right),\hspace{1cm} J_3:=\left(
\begin{array}{cc}
  0 & 1 \\
  -1 & 0 \\
\end{array}\right).
$$

The one-parameter group generated by $J_3$ is compact, while those
 generated by $J_1$ and $J_2$ are non-compact. Let $\lambda>0$. The
 manifold
$$N:=\{\xi_1 J^1+\xi_2 J^2+\xi_3 J^3;\xi_3^2-\xi_1^2-\xi_2^2=\lambda^2, \xi_3>0\}$$
is a coadjoint orbit. Then the generalized function
$\vol_N(\phi_1J_1+\phi_2J_2+\phi_3J_3)$ is given by an invariant
locally $L_1$-function, analytic outside
$\phi_1^2+\phi_2^2-\phi_3^3=0.$

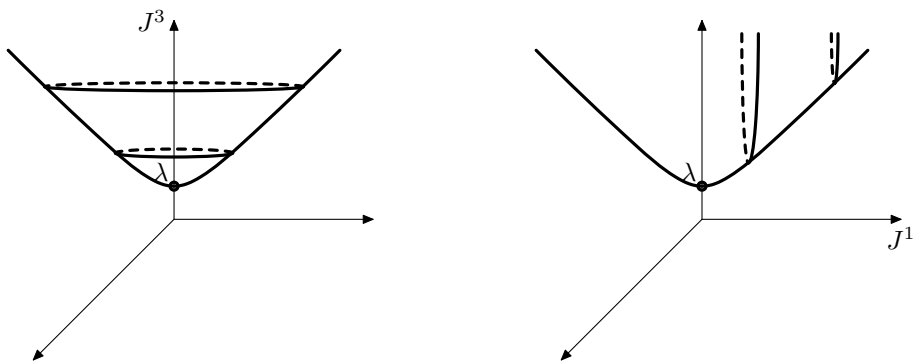
\vspace{1cm}

\includegraphics{orbit1.01} \hspace{2cm}\includegraphics{orbit1.02}

$$\vol_N(\phi_3 J_3)=-\frac{e^{i\lambda\phi_3}}{2i\phi_3},
\hspace{4cm}\vol_N(\phi_1 J_1)= \frac{e^{-|\lambda
\phi_1|}}{2|\phi_1|}.$$

 The  formula for the
generator $J_3$ of a compact group action is in agreement with the
``fixed point formula philosophy". 
The  formula  for $J_1$ is difficult to explain within  a general
framework. Indeed, the non-compact group $\exp(\phi_1 J_1)$
 acts freely on $N$; however, the value of the function
$\vol_N(\phi_1 J_1)$ is non-zero even though there are no fixed
points on $N$. In \cite{lib04}, $N$ is embedded in the cotangent
bundle of the  Riemann sphere $M:=P_1(\C)$, and a subtle argument
of deformation to fixed points of $J_1$ in $M$ ``explains" the
formula for $\vol_N(\phi_1 J_1)$.

\subsection{Equivariant cohomology groups}\label{equicoho}
After having defined $\CH^{\infty}(\g,N)$, I will move on to the
definition of  two other equivariant cohomology groups.

{\bf Cartan's complex.}  Here we consider, for a $K$-manifold $N$,
the space  $\CA^{\rm pol}(\k,N):=(S(\k^*)\otimes \CA(N))^K$
 of equivariant forms $\form(\phi)$
depending polynomially on $\phi$. The corresponding cohomology
space $\CH^{\rm pol}(\k,N)$ is a $\Z$-graded algebra, where
elements of $\k^*$ have degree two, and differential forms their
exterior degree. If $N$ is a vector space with linear action of
$K$, then $\CH^{\rm pol}(\k,N)=S(\k^*)^K$. A basic theorem of H.
Cartan says: if $K$ acts on a compact manifold $M$ with finite
stabilizers, then $\CH^{\rm pol}(\k,M)=\CH^*(M/K)$.

If $N$ is non-compact, we can also consider
 the space  $\CA^{\rm pol,cpt}(\k,N):=(S(\k^*)\otimes \CA^{\rm cpt}(N))^K$
  of equivariant forms $\form(\phi)$
 which are  compactly supported on $N$. We denote by
$\CH^{\rm pol,cpt}(\k,N)$ the corresponding cohomology space.
Integration is well defined on it if $N$ is oriented and  the
result of integration $\int_N\form(\phi)$ is a polynomial function
on $\k$, invariant under the adjoint action of $K$ on $\k$.

If $N$ is a vector space, there exists a unique element ${\rm
Thom}(\phi)\in \CH^{\rm pol,cpt}(\k,N)$ with integral equals to
$1$.
\bigskip
Let us give the formula for $N:=\R^2$ with action of $S^1$.
$\bullet$ $N:=\R^2$. Let $\chi$ be any smooth compactly supported
function on $\R$ such that $\chi(0)=1$. Then
\begin{equation}\label{Thom}
{\rm
Thom}_{\chi}(\phi):=\frac{-1}{2\pi}(\phi\,\chi(x^2+y^2)+2\chi'(x^2+y^2)
dx\wedge dy)
\end{equation}
 is a representative of ${\rm Thom}(\phi)$.

If $N$ is a vector space, a representative of ${\rm Thom}(\phi)$
with ``Gaussian look" is given by Mathai-Quillen in
\cite{mat-qui}.

Details  on Cartan's theory  and  further developments can be
found in  the stern monograph (which contains treasures)
\cite{duf-kum-ver}, or in  the attractive  book \cite{gui-ste99}.
This de Rham point of view for topological equivariant cohomology
seems to be adapted only to smooth spaces. However, the use of
equivariant Poincar{\'e} dual allows us  to  work  on algebraic
varieties, where the Joseph polynomials and the Rossmann
localization formula (see \cite{ros2}) are important tools. For
lack of space, I will not pursue this topic. Let me also mention
the theory  of
 equivariant Chow groups for algebraic actions on algebraic
varieties defined over any field, initiated by Totaro
 and developed by Edidin-Graham   and Brion.

{\bf Generalized coefficients}\,\,\cite{duf-kum-ver}. An
equivariant  form $\alpha(\phi)$ with $C^{-\infty}$ coefficients
is a generalized function on $\g$ with values in $\CA(N)$. Thus
for any  smooth function $F$  on $\g$ with compact support, the
integral $\int_{\g}\form(\phi)F(\phi) d\phi$ is a differential
form  on $N$. We denote by $\CA^{-\infty}(\g,N)$ the space of such
forms. If $N:=\bullet$ is a point, an equivariant form with
$C^{-\infty}$ coefficients is just an element of
$(C^{-\infty}(\g))^G$, that is, an invariant generalized function
on $\g$. The operator $D$ is well defined on
$\CA^{-\infty}(\g,N)$, and  we denote the corresponding cohomology
space by $\CH^{-\infty}(\g,N)$. It is a module over
$\CH^{\infty}(\g,N)$. If $K$ acts freely on $N$, the natural image
of $\CH^{\infty}(\k,N)$ in $\CH^{-\infty}(\k,N)$ is equal to $0$.

\begin{example}
 Let $M:=S^1=\{e^{i\theta}\}$. The group  $S^1:=\{e^{i\phi}\}$ acts freely on $M$ by rotations. Let $\g:=\R J$ be the Lie algebra of $S^1$. Then $\CH^{-\infty}(\g,M)=\C v$,
where $\phi v=0$. A representative of $v$, still denoted by $v$,
is the closed equivariant form
$$v(\phi):=\Dirac(\phi)d\theta.$$
 Note that  $\int_M v(\phi)=(2\pi)\Dirac(\phi).$

On the other hand, we have  $1=-i\,D(Y^+(\phi)d\theta)$ so that
$$1=0 \hspace{1cm} {\rm in}\,\, \CH^{-\infty}(\g,M).$$
Thus the image of $\CH^{\infty}(\g,M)=\C \cdot 1$ in
$\CH^{-\infty}(\g,M)$ vanishes.
\end{example}

\subsection{Localization or $1=0$}\label{wit}
 Let $N$ be a $K$-manifold and let
 $\kappa$ be a $K$-invariant vector field, tangent
 to  $K$-orbits. If $\nu:N\to \k$ is a $K$-invariant map, then $\kappa$, defined by
\begin{equation}\label{kappa}
  \kappa_m:=\frac{d}{d \epsilon}\exp(\epsilon \nu(m)).m|_{\epsilon=0}
  \end{equation}
  is such a vector field.

   Let $C$  be the set of zeroes of $\kappa$. Via a $K$-invariant Riemannian structure $( \bullet,\bullet) $ on $N$, identify $\kappa$
 with
 the
 $K$-invariant $1$-form on $N$: $\langle\kappa,\bullet\rangle:=(\kappa,\bullet).$
 Witten considers the exact
equivariant form $D\kappa(\phi)=-\langle \kappa,V\phi\rangle +
d\kappa.$ From our tangential hypothesis, $\phi\mapsto \langle
\kappa_m,V\phi\rangle $ is a non-zero element of $\k^*$ when $m$
is not in $C$.

Let $\form(\phi)\in \CH^{\infty}(\k,N)$, compactly supported on
$N$. For any test function $F(\phi)$ on $\k$ and any $a$ in $\R$,
we have the equality
\begin{equation}\label{defwitten}
\int\int_{N\times \k} \form(\phi)F(\phi)d\phi=\int_N \int_{\k}
e^{-ia D\kappa(\phi)}\form(\phi) F(\phi)d\phi.\end{equation}

When $a$ tends to infinity, standard estimates on Fourier
transforms shows that the differential form  $\int_{\k} e^{-ia
D\kappa(\phi)}\form(\phi) F(\phi)d\phi$ becomes very small outside
$C$.

 Inspired by Witten's deformation argument, Paradan proves that outside $C$, the constant $1$ is equal to $0$ in $\CH^{-\infty}(\k, N-C)$.
\begin{theorem}\cite{par00}\label{par}(Paradan)
On $N-C$, the integral
$$B(\phi):=i\int_{0}^{\infty}e^{-ia D\kappa(\phi)}\kappa da$$ is a well defined  element of
$\CA^{-\infty}(\k, N-C)$ and we have $1=D(B(\phi))$. Thus
$$1=0\,\,\, {\rm  in}\,\, \CH^{-\infty}(\k, N-C).$$
\end{theorem}
Indeed, intuitively $B(\phi)=\frac{\kappa}{D\kappa(\phi)}$, so
that $DB(\phi)= \frac{D\kappa(\phi)}{D\kappa(\phi)}=1$.

Multiplying an element $\alpha(\phi)$ of $\CH^{\infty}(\k,N)$  by
$1$, we see that $\alpha(\phi)$ vanishes on $N-C$. In the next
proposition, we give an explicit representative of $\alpha$ with
support near $C$.
\begin{proposition}\cite{par00}\label{parbis}
 Let $\chi$ be a $K$-invariant function on $N$ supported on a
small neighborhood of $C$ and such that $\chi=1$ on a smaller
neighborhood of $C$. Let
$$P(\phi):=\chi+ d\chi\wedge B(\phi).$$
Then $P(\phi)$ is a closed equivariant form in
$\CA^{-\infty}(\k,N)$ supported near  $C$. Furthermore, we have
the equation  in $\CA^{-\infty}(\k,N)$:
$$ P=1+D((\chi-1)B).$$
Thus, if  $\alpha(\phi)\in \CH^{\infty}(\k,N)$, then
$P(\phi)\alpha(\phi)$ is supported near $C$ and equal to
$\alpha(\phi)$ in $\CH^{-\infty}(\k,N)$.
\end{proposition}

In the basic examples $\R^2$ or $T^*S^1$ with action of $S^1$, and
$\kappa$ appropriately chosen, the forms $B(\phi)$  and $P(\phi)$
are easy to calculate.

$\bullet$ Consider $N:=\R^2$ with
$\kappa:=y\partial_x-x\partial_y$. On $\R^2-\{[0,0]\}$, in polar
coordinates $r,\theta$, we compute that $B(\phi)=-iY^{+}(\phi)
d\theta.$ Thus, if $\chi$ is a smooth function with compact
support on $\R$ and equal to $1$ in a neighborhood of $0$, then
 $$P(\phi)=(2i\pi)Y^{+}(\phi){\rm Thom}_\chi(\phi),$$
 where ${\rm Thom}_\chi(\phi)$ is defined by Formula (\ref{Thom}).
Note that the integral of $P(\phi)$ on $N$ is $(2i\pi)Y^+(\phi)$.

$\bullet$ Consider $N:=T^*S^1$ with $\kappa:=-t \partial_\theta$.
Then in coordinates $t,\theta$,
\begin{eqnarray*}\label{parphi}
B(\phi)&=&-iY^{+}(\phi) d\theta \hspace{1cm} {\rm if}\,\, t>0,\\
B(\phi)&=&iY^{-}(\phi) d\theta  \hspace{1cm} {\rm if}\,\, t<0.
\end{eqnarray*}
If $\chi$ is a smooth function with compact support on $\R$ and
equal to $1$ in a neighborhood of $0$, then
$$P(\phi)=\chi(t)+  \chi'(t) dt\wedge B(\phi).$$ Note that the
integral of $P(\phi)$ on $N$ is
$(2i\pi)(Y^+(\phi)+Y^-(\phi))=(2i\pi)\Dirac(\phi)$.

\bigskip

For the sake of simplicity, assume that $N$ is compact. Consider
the form $P\in \CH^{-\infty}(\k,N)$ constructed in Proposition
\ref{parbis} and supported near  the set $C$ of zeroes of
$\kappa$. We write $C=\cup C_F$ where $C_F$ are the connected
components of the set $C$. Write $P=\sum_F P_F$ where $P_F$ is
compactly supported on a small neighborhood $U_F$ of $C_F$.
Proposition \ref{parter} reduces the calculation of the integral
of $\alpha(\phi)$ on $N$ to calculations near $C$. We obtain the
following  localization theorem.

\begin{theorem}\cite{par00}\label{parter}. Consider  an equivariant class $\alpha(\phi)\in
\CH^{\infty}(\k,N)$. For any component $C_F$ of the set $C$, let
$\alpha_F(\phi)\in \CH^{\infty}(\k, U_F)$ an equivariant class
 equal to $\alpha(\phi)$ on $U_F$. Then
 $$\int_N
 \alpha(\phi)=\sum_{C_F}\int_{U_F}P_F(\phi)\alpha_F(\phi).$$
\end{theorem}
In this localization theorem, each local contribution
$\int_{U_F}P_F(\phi)\alpha_F(\phi)$ is a generalized function on
$\k^*$. Thus the Fourier transform of each local contribution has
a meaning, under a moderate growth condition for $\alpha$.

 As an application, we recover the
exact stationary phase, and more generally   the ``abelian"
localization formula, with the following tool. For a $S^1$-action
with generator $J$, we choose $\kappa:= J$, so that $C$ is the set
of fixed points of the one parameter group  $\exp(\phi J)$. We
obtain the following result that we state in the case of isolated
fixed points.
\begin{theorem}\label{BVAB}
(\cite{ber-ver82-1}, \cite{wit82}, \cite{ati-bot84}) Let $S^1$
acting  on a compact manifold  $M$ with isolated fixed points. Let
$\form(\phi)$ be a closed equivariant form with $C^{\infty}$
coefficients. Then
$$(2\pi)^{-\frac{\dim M}{2}}\int_M \form(\phi)
=\sum_{p\in \{\rm fixed\, points\}}
\frac{i_p^*\form(\phi)}{\sqrt{\det_{T_pM}L_p(\phi)}}.$$
\end{theorem}

\section{Applications and Conjectures}\label{appli}

\subsection{Integrals on reduced spaces}
\subsubsection{Reduced spaces}\label{volred}

 Let $N$ be a Hamiltonian $K$-manifold. Assume that $\xi\in \k^*$
 is a regular value of the moment map $\mu$ and
let $K_\xi$ be the stabilizer of $\xi$. Then $K_\xi$ acts with
finite stabilizers in ${\bf \mu}^{-1}(\xi)$ so that ${\bf
\mu}^{-1}(\xi)/K_\xi$ is a symplectic orbifold, called the reduced
space at $\xi$ and denoted by $N_\xi$. We denote by $s_\xi$ the
number of elements of the stabilizer  of a generic point in
$\mu^{-1}(\xi)$. If $\xi=0$, we also denote $N_0=\mu^{-1}(0)/K$ by
$N//K$. When $N$ is a projective manifold, then $N//K$ is the
quotient in the sense of Mumford's geometric invariant theory (see
chapter 8.2 \cite{mum-fog-kir}). By considering the symplectic
manifold $N\times (K\cdot(-\xi))$ (the shifting trick), we may
always consider reduction at $0$.

If $0$ is a regular value, Kirwan associates to an equivariant
closed form $\alpha(\phi)$ on $N$ a cohomology class $\alpha_{\rm
red}$ on $N//K$:  $\form(\phi)|_{{\bf \mu}^{-1}(0)}$ is equivalent
to the pull-back of $\alpha_{\rm red}$. The Kirwan map
$\chi:\CH^*_K(N)\to \CH^*(N//K)$ is surjective, at least when $N$
is compact.

The following result relates the equivariant volume of $M$ to
volumes of reduced spaces.
\begin{proposition}\cite{dui-hec}\label{dudu}
If $M$ is a $K$-Hamiltonian manifold, then
$$\vol_M(\phi)=\int_{\k^*} e^{i\langle \xi,\phi\rangle }\vol(M_\xi)d\xi.$$
\end{proposition}

This theorem holds also if $N$ is a $K$-Hamiltonian manifold with
proper moment map, under some convergence conditions. As shown by
Formula (\ref{2SD}) (Section \ref{equivol}), if a torus $T$ acts
on a vector space $N$ with weights $\beta_a\in \t^*$, all
contained in a half-space, then the equivariant volume
$\vol_N(\phi)$ is the boundary value of
$\frac{1}{\prod_{a}(-i\beta_a(\phi))}$. Its Fourier transform is
the
  convolution $H$
of the Heaviside distributions supported on the half-lines
$\R^+\beta_a$. Computing  volumes of the reduced manifolds $N_\xi$
is the
 same as computing the value of $H$ at a point $\xi\in \t^*$. In
 Section \ref{Heaviside}, we will explain how to do it using
iterated residues.

 In the next
section, we explain Witten's generalization of Proposition
\ref{dudu}.

\subsubsection{Witten's localization theorem}
Assume that $M$ is a compact $K$-Hamil-tonian manifold with moment
map ${\bf \mu}: M\to \k^*$. We choose a $K$-invariant
identification $\k^*\to \k$ given by a $K$-invariant inner
product. The vector field $\kappa$ defined by
$\kappa_m:=\frac{d}{d\epsilon}\exp(-\epsilon {\bf \mu}(m))\cdot
m|_{\epsilon=0}$ is $K$-invariant. We refer to this particularly
important vector field as the Kirwan vector field. In this case,
the set $C$ of zeroes of $\kappa$ is the set of critical points of
the invariant function $\|{\bf \mu}\|^2$ on $M$. One connected
component of $C$ is the set ${\bf \mu}^{-1}(0)$ of zeroes of the
moment map (if not empty). The following theorem follows from
Witten's deformation argument: Formula (\ref{defwitten}) in
Section \ref{wit}.

\begin{theorem}\label{Witten}(Witten \cite{wit92})
Let $M$ be a compact  Hamiltonian $K$-manifold and  $p(\phi)$ an
equivariantly closed form with polynomial coefficients. Assume
that  $0$ is  a regular value of the moment map. Then
$$\int_{ \k} \left(\int_Me^{i\Omega(\phi)}p(\phi)\right)d\phi=
s_0(2i\pi)^{\dim \k} \vol(K)\int_{M//K}e^{i\Omega_{\rm red}}p_{\rm
red}.$$

\end{theorem}

Let me explain the meaning of the integral on the left. Let
$I_M(\phi):=
   \int_M e^{i\Omega(\phi)}p(\phi)$. This is an analytic function
on $\k$ with at most polynomial growth. We compute $\int_{ \k}
e^{i \langle \xi,\phi \rangle}I_M(\phi)d\phi$ in the sense of
Fourier transform. This Fourier transform is a polynomial near
$\xi=0$ (this is part of the theorem). The  left-hand side
$\int_{\k}I_M(\phi)d\phi$ is by definition the value of this
polynomial at $\xi=0$.

The theorem above is used to compute integrals on reduced spaces.
Indeed, the right hand side of the equality is the integral of a
cohomology class over the reduced space $M//K$ of $M$, which is
difficult to compute. Instead, we first compute an equivariant
integral on the original space $M$ (easy to do thanks to  the
usual reduction to the maximal torus $T$ and the abelian
localization formula). Then we have to compute the value of the
Fourier transform  of $I_M(\phi)$ at the point $0$. This in turn
demands  the computation of the value of the convolution of
Heaviside distributions at some explicit points of $\t^*$: the
images
by $\mu$ of the fixed points of the action of $T$ on $M$.

Using different methods, other proofs and refinements  to Witten's
theorem have been given
 (\cite{jef-kir95}, \cite{ver96-1}, \cite{par00}, \cite{sawin}).
Let us recall  Paradan's method. We apply  Theorem \ref{parter} to
the form $\alpha(\phi)=e^{i\Omega(\phi)}p(\phi)$. Here $C_F$
varies over the connected components of the set of critical points
of $\|\mu\|^2$. The image of a connected component $C_F$ by the
moment map $\mu$ is a $K$-orbit $K\beta$. The set
$C_0:=\mu^{-1}(0)$ projecting on $0$ is one connected component of
$C$ (if non-empty). The Fourier transform of $\int_M
e^{i\Omega(\phi)}p(\phi)P_F(\phi)$  when $C_F$ projects on
$K\beta$ with $\beta\neq 0$ is supported on $\|\xi\|\geq
\|\beta\|$. Thus the value of the Fourier transform of $\int_M
e^{i\Omega(\phi)}p(\phi)$ at $0$ comes only from $\int_M
e^{i\Omega(\phi)}p(\phi)P_0(\phi)$ and requires only local
knowledge of our data near $\mu^{-1}(0)$. To summarize, in
Witten's localization formula, the Fourier transform of the local
terms arising from components different from $C_0$ are moved away
from our focus of attention: the point $0$ in $\k^*$.

 These facts are illustrated
in the example below. This also shows   that local calculations
near critical points  essentially  reduce to $\R^2$ or $T^*S^1$.

\begin{figure}[!h]
\begin{center}
\includegraphics{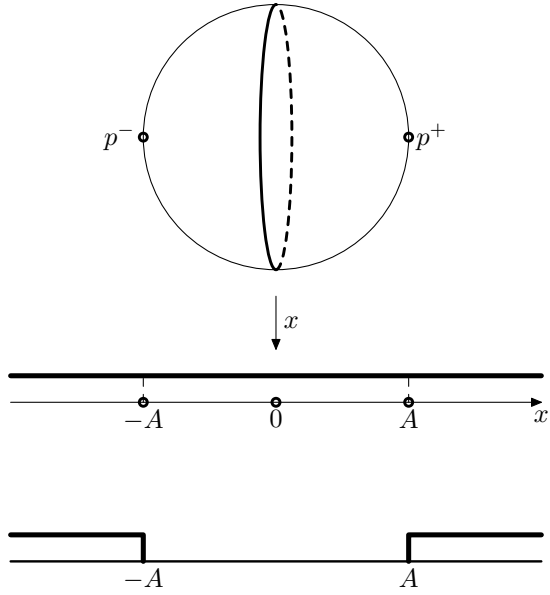}
\end{center}
\caption{Decomposition of equivariant volumes}\label{witint}
\end{figure}

\begin{example}\label{witonsphere}
Return to Example \ref{sphere} of the sphere
$M:=\{x^2+y^2+z^2=A^2\}$, with
 moment map ${\bf \mu}(x,y,z)=x$.
The  critical values of $x^2$ are $0,A,-A$. The set of critical
points has three connected components: the circle $C_0$ drawn in
black in Figure \ref{witint}, $\{p^+\}$ and $\{p^-\}$. The normal
bundle to $C_0$ is identified with $T^*S^1$ and the normal bundles
to $p^+,p^-$ with  $\R^2$. Let
$$\vol_M(\phi)=\frac{1}{2i\pi}\int_Me^{i\Omega(\phi)}.$$
Using the Kirwan vector field, we  obtain a decomposition
$\vol_M(\phi)=v_0(\phi)+v_{p^-}(\phi)+v_{p^+}(\phi)$ with
\begin{eqnarray*}
v_0(\phi):&=&\frac{1}{2i\pi}\int_M
e^{i\Omega(\phi)}P_0(\phi)=\delta_0(\phi),\\
v_{p^-}(\phi):&=&\frac{1}{2i\pi}\int_M
e^{i\Omega(\phi)}P_{p^-}(\phi)=-e^{-i\phi A}Y^{-}(\phi),\\
v_{p^+}(\phi):&=&\frac{1}{2i\pi}\int_M
e^{i\Omega(\phi)}P_{p^+}(\phi)=-e^{i\phi A}Y^{+}(\phi).
\end{eqnarray*}
 This decomposition corresponds to the cone
decomposition of the interval $[-A,A]$ described in  Figure
\ref{decomint} in Section \ref{geo}.
\end{example}

\subsection{Index of transversally elliptic
operators}\label{transver}

 Consider a compact even-dimensional oriented manifold $M$.
For the sake of simplicity, we assume $M$ provided with an almost
complex structure. We choose an Hermitian metric $\|\xi\|^2$ on
$T^*M$. For $[x,\xi]\in T^*M$, the symbol of the Dolbeault-Dirac
operator $\overline \partial +\overline \partial^*$ is the
Clifford multiplication  $c(\xi)$ on the complex vector bundle
$\Lambda T_x^*M$. It is invertible for $\xi\neq 0$, since
$c(\xi)^2=-\|\xi\|^2$. Let $\CE$ be an auxiliary vector bundle
over $M$, then $c_{\CE}([x,\xi]):=c(\xi)\otimes {\rm Id}_{\CE_x}$
defines an element of the ${\mathbf K}$-group of $T^*M$. Assume
that a compact group $K$ acts on $M$ and $\CE$. Now the
topological index ${\rm Index}(c_{\CE})$ of $c_{\CE}\in {\mathbf
K}_K(T^*M)$ is an invariant function on $K$ (which computes the
equivariant index of the $K$-invariant operator $\overline
\partial_{\CE} +\overline \partial^*_\CE$).
The index theorem of  Atiyah-Segal-Singer expresses ${\rm
Index}(c_{\CE})(k)$ ($k\in K)$ in terms of the fixed  points of
$k$ on $M$. We constructed (see \cite{ber-get-ver})  the
equivariant Chern character ${\rm ch}(\phi,\CE)$ of the vector
bundle $\CE$ and the equivariant Todd class ${\rm Todd}(\phi,M)$
such that (for $\phi$ small)
\begin{equation}\label{index}
{\rm Index}(c_{\CE})(\exp \phi)=(2i\pi)^{-(\dim M)/2}\int_M {\rm
ch}(\phi, \CE) {\rm Todd}(\phi,M).
\end{equation}

For $\phi=0$, this is the Atiyah-Singer formula. Formula
(\ref{index}) is a ``delocalization" of the Atiyah-Segal-Singer
formula. The delocalized index formula (\ref{index}) can be
adapted to new cases such as:
\begin{itemize}
\item Index of transversally elliptic operators.

\item $L^2$-index of some elliptic operators on some non-compact
manifolds (as in Narasimhan-Okamoto, Parthasarathy, Atiyah-Schmid,
Connes-Moscovici).
\end{itemize}

Indeed,  in these two contexts, the index exists in the sense of
generalized functions but cannot be always computed in terms of
fixed  point formulae.

 Recall  Atiyah-Singer's definition of transversally elliptic operators
 (see \cite{ati74}).
Let $N$ be a  $K$-manifold and  $T^*_K N$ be  the conormal bundle
to $K$-orbits. A  transversally elliptic pseudo-differential
operator $S$  is elliptic in the directions normal to the
$K$-orbits. Thus $S$ together with the action of the Casimir of
$\k$ defines an elliptic system, and the space of solutions of $S$
decomposes as a Hilbert direct sum of finite-dimensional spaces of
$K$-finite solutions. The symbol of $S$ defines an element
  $\sigma(S)$ of ${\mathbf K}_K(T^*_KN)$. The index of the operator $S$ is the character
of $K$ in the virtual vector space obtained as difference  of
$K$-finite solutions of $S$ and its adjoint. This is an invariant
generalized function on $K$. In \cite{ber-ver96-2}, we gave a
cohomological formula for the index of   $S$ in terms   of
$\sigma(S)\in {\mathbf K}_K(T^*_KN)$, as an equivariant integral
on $T^*N$ in the spirit of the delocalized formula (\ref{index}).
This result was inspired by Bismut's ideas on delocalizations
\cite{bis86} and Quillen's superconnection formalism.

The following example shows that, contrary  to the melancholy
remark of Atiyah about his work on transversally elliptic
operators (page 6, vol 4, \cite{ati0}), there are many
transversally elliptic bundle maps of great interest.

 Consider a
$K$-manifold $N$ with  a $K$-invariant vector field $\kappa$
tangent to orbits. As before, we assume that $N$ is provided with
a $K$-invariant almost complex structure and Hermitian metric. We
still denote by $c(\xi)$ the Clifford action of $\xi\in T^*_xN$ on
the complex space $\Lambda T_x^*N$. The analogue in $K$-theory of
Witten's deformation is the bundle map

\begin{equation}\label{paroperator}
c_{\kappa,\CE}([x,\xi]):=c(\xi-\kappa_x)\otimes {\rm Id}_{\CE_x},
\end{equation}
  defined by Paradan \cite{par2}.
Note that $c_{\kappa,\CE}([x,\xi])$ is invertible except if $
\xi=\kappa_x$. If furthermore $[x,\xi]\in T^*_K N$, then  $\xi=0$
and $\kappa_x=0$. Indeed, by our hypothesis, under  identification
of $T^*N$ with $TN$, $\kappa_x$ is  tangent to $Kx$ while $\xi$ is
normal to $Kx$.

 When $N$ is compact, $c_{\kappa,\CE}$ is transversally elliptic and
 equal in
${\mathbf K}
$-theory to  the elliptic symbol $c_{\CE}$, via the
deformation
 $c(\xi-a\kappa_x)\otimes 1_\CE,$ for $a\in [0,1]$.
Under the conditions stated below, Paradan's construction
  defines a transversally elliptic element
 even if $N$ is a non-compact manifold. See also the construction by M. Braverman \cite{bra}
 of a related operator.
\begin{proposition}\cite{par2}\label{Dolbeault}
 Assume that the set $C$ of zeroes of $\kappa$ is compact. Then $c_{\kappa,\CE}$ is
transversally elliptic on $T^*N$ with support the zero section
$[C,0]$.
\end{proposition}

Recall the closed equivariant form $P$  on $N$  supported on a
neighborhood of $C$ constructed with the help of $\kappa$ in
Proposition \ref{parbis}. Then

\begin{theorem}\cite{par-ver}
Near the identity element $1$ of $K$, the index of
$c_{\kappa,\CE}$ is given by the formula
\begin{equation}\label{conjectureOTE}
{\rm Index}(c_{\kappa,\CE})(\exp \phi)= (2i\pi)^{-(\dim N)/2}
\int_{N} {\rm ch}(\phi,\CE) {\rm Todd}(\phi,N)
P(\phi)\end{equation} and by similar integral formulae over $N_s$
near any point $s\in K$.
\end{theorem}

When $M$ is compact, Formula (\ref{conjectureOTE})  reduces to
Formula (\ref{index})  since $P(\phi)$ is equal to $1$ in
cohomology. But even in this case, Formula (\ref{conjectureOTE})
has important implications, as the symbol $c_{\CE}$ is broken into
several parts according to the connected components of $C$:
$c_{\CE}=\sum_{F} c_{\CE,F}$  where $c_{\CE,F}$ is supported on
$[C_F,0]$. Thus
$${\rm Index}(c_\CE)=\sum_{F} {\rm Index}(c_{\CE,F}).$$
Each local contribution ${\rm Index}(c_{\CE,F})$ is well defined
as a  character of an infinite-dimensional representation of $K$.
This was one of the  motivations of Atiyah and Singer  for
introducing transversally elliptic operators.

 As in the Witten localization formula, this allows in
important cases to compute the invariant part ${\rm
Index}(c_\CE)^K$ through considering only the contribution of
$C_0$. The  Fourier series attached to the other components do not
interfere with  our focus of attention: the multiplicity of the
trivial representation. This fact is illustrated in the example
below.

\begin{figure}[!h]
\begin{center}
\includegraphics{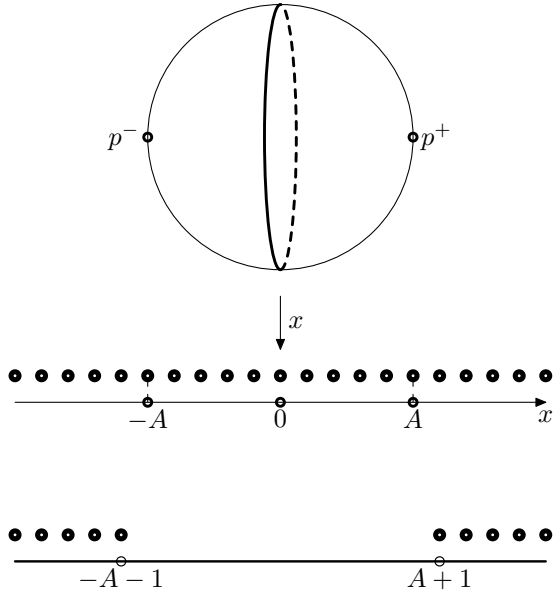}
\end{center}
\caption{Decomposition of equivariant indices}\label{witindices}
\end{figure}

\begin{example}\label{ellipticonsphere}
Return to Example \ref{witonsphere}. Let $A$ be a positive
integer. We identify $P_1(\C)$ with $M_A:=\{x^2+y^2+z^2=A^2\}$
through the map $$[z_1,z_2]\mapsto
(A\frac{|z_1|^2-|z_2|^2}{|z_1|^2+|z_2|^2},2 A\frac{\Re(z_1
\overline{z_2})}{|z_1|^2+|z_2|^2},2 A\frac{\Im(z_1
\overline{z_2})}{|z_1|^2+|z_2|^2}),$$ the action $(e^{i
\phi}z_1,z_2)$ becoming the rotation around the  $x$-axis. We
consider the  Dolbeault-Dirac  operator $D_{2A}$ on $P_1(\C)$ with
solution space $\oplus_{j+k=2A} \C z_1^j z_2^k$. Twisting the
action by $e^{i \phi A}$, its equivariant index is $\sum_{k=-A}^A
q^k$ with $q:=e^{i\phi}$. Using the Kirwan vector field, we
decompose $D_{2A}=D_0+D_{p^+}+ D_{p^-}$ into the sum of three
transversally elliptic operators with support $[C_0,0]$,
$[p^+,0]$, $[p^-,0]$, respectively. To compute the index of $D_0$,
we are led to compute the set of solutions of the Dolbeault
operator on the complex manifold $\C/\Z=S^1\times \R$, the action
of $S^1=\R/2\pi\Z$ being by translations, and we obtain all
functions $e^{i k z}$ for any $k\in \Z$. Thus $${\rm Index
}(D_0)=\sum_{k=-\infty}^{\infty}q^k.$$ Near the fixed points
$p^+,p^-$, we obtain  the index of the lift of the operators
$\overline\partial^{\pm}$ (see \cite{ati74}) on $\C$ (shifted):
$${\rm Index}
(D_{p^+})=-\sum_{k= A+1}^{\infty}q^k,\hspace{1cm} {\rm Index
}(D_{p^-})=-\sum_{k= -\infty}^{-A-1}q^k.$$

The equality  $${\rm Index} (D_{2A})={\rm Index}( D_{0})+{\rm
Index}( D_{p^+})+{\rm Index}( D_{p^-})$$ is   Formula
(\ref{geotrans}) in Section \ref{geo}.
\end{example}

It might happen that the integral $\int_N {\rm ch}(\phi,\CE) {\rm
Todd}(\phi,N)$ over our non-compact manifold $N$ is already
convergent in the distributional sense, and as $P=1$ in
cohomology, it might happen, modulo the convergence of the
boundary term, that the following equality holds
$${\rm Index}(c_{\kappa,\CE})(\exp \phi)=  (2i\pi)^{-(\dim N)/2}
\int_N {\rm ch}(\phi,\CE) {\rm Todd}(\phi,N).$$

This is indeed the case for discrete series. To state the result,
we rephrase the preceding constructions in the spin context. If
$N$ is an even-dimensional oriented spin manifold, and $\CE$ a
twisting vector bundle, we denote by $\sigma(\xi)$ the Clifford
action of $\xi\in T_x^*N$ on spinors, and by $\sigma_\CE$ the
symbol of the twisted Dirac operator $D_\CE$. If $M$ is a compact
$K$-manifold, the equivariant index of $D_\CE$ is given by a
formula similar to (\ref{index}) :
\begin{equation}\label{indexdirac}
{\rm Index}(\sigma_{\CE})(\exp \phi)=(2i\pi)^{-(\dim M)/2}\int_M
{\rm ch}(\phi, \CE) {\rm \hat A}(\phi,M),
\end{equation}
where the equivariant class $\hat A$  replaces the equivariant
Todd class.

Under the same hypothesis as in Proposition \ref{Dolbeault}, the
bundle map
$$\sigma_{\kappa,\CE}([x,\xi])=\sigma(\xi-\kappa_x)\otimes {\rm
I_{\CE_x}}$$ is transversally elliptic and its equivariant index
is a generalized function on $K$.

Let $G$ be a real reductive Lie group with maximal compact
subgroup $K$. We assume that the maximal torus $T$ of $K$ is a
maximal torus in $G$. Let $N:=G\lambda$ be the orbit of a regular
admissible element $\lambda\in\t^*$. Harish-Chandra associates to
$\lambda$ a representation of $G$, realized as the $L^2$-index of
the twisted Dirac operator $D_\lambda$. The moment map $\mu$ for
the $K$-action on $N$ is the projection $G\lambda\to \k^*$ and the
set $C$ of zeroes of the Kirwan vector field $\kappa$ is  easy to
compute in this case: it consists of the compact orbit $K\cdot
\lambda$.

\begin{theorem}\label{discrete}(Paradan \cite{par3})
The character of the discrete series $\Theta^G(\lambda)$
restricted to $K$ is the index of the transversally elliptic
element $\sigma_{\kappa,\CL_\lambda}$ on $N$.
\end{theorem}
Here $\CL_\lambda$ is the Kostant line bundle
$G\times_{G(\lambda)}\C_{\lambda}$ on $N=G/G(\lambda)$. A
calculation of the index of $\sigma_{\kappa,\CL_\lambda}$ (which
is supported
  on $K\cdot \lambda$) leads immediately to  Blattner's
formula for $\Theta^G(\lambda)|_K$.

\subsection{Quantization and symplectic quotients}\label{QR} Let
$N$ be a  $G$-manifold  ($N,G$ non-necessarily compact), and $\CE$
a $G$-equivariant vector bundle on $N$ with $G$-invariant
connection $\nabla$. We can then construct the closed equivariant
form ${\rm ch}(\phi,\CE)$ (\cite{ber-ver82}, \cite{bot-tu}). For
the sake of simplicity, I assume the existence of a $G$-invariant
complex structure. Then I conjectured (under additional conditions
that I do not know how to formulate exactly, see attempts  in
\cite{ver92})

{\bf Conjecture:} There exists a
 representation $Q(N,\CE)$ of $G$ such that the character
$\Tr_{Q(N,\CE)}(g)$ is given by the formula
\begin{equation}\label{ecm}
\Tr_{Q(N,\CE)}(\exp \phi)= (2i\pi)^{-(\dim N)/2} \int_N {\rm
ch}(\phi,\CE) {\rm Todd} (\phi,N) \end{equation}
 near $1\in G$ and by a similar integral  formula  over $N_s$ near
any elliptic point $s$ of $G$.

Thus, via integration of equivariant cohomology classes, it should
be possible to define a push-forward map from a generalized
${\mathbf K}$-theory of vector bundles with connections on
$G$-manifolds  to invariant generalized functions on $G$, under
some convergence conditions, and assuming the existence of a
suitable equivariant Todd class.

\begin{remark}\label{orbit}
 When $N$ is a coadjoint admissible regular orbit of
any real algebraic Lie group $G$ and $\CE$ the Kostant half-line
bundle, Formula (\ref{ecm}), with the $\hat A$ class instead of
the Todd class, becomes  Kirillov's universal formula \cite{kiri5}
for characters (proved by Kirillov for compact and nilpotent
groups, by Duflo, Rossmann, Bouaziz, Khalgui, Vergne,... for any
real algebraic group). If $N,G$ are compact, Formula (\ref{ecm}),
with $\hat A$ instead of ${\rm Todd}$, is the equivariant index
formula for the Dirac operator twisted by $\CE$. Thus Formula
(\ref{ecm}), modified as in \cite{ver92}, is a fusion of the
Kirillov  universal character formula and of the formulae of
Atiyah-Segal-Singer for indices of twisted Dirac operators.
\end{remark}

Now let  $(M,\Omega)$ be a compact symplectic manifold with
Hamiltonian action of  a compact group $K$. We assume the
existence of a $K$-equivariant line bundle $\CL$  on $M$ with
connection $\nabla$ of curvature equal to $i\Omega$. In other
words, $M$ is prequantizable in the sense of \cite{kos} and we
call $\CL$ the Kostant line bundle. We take an almost complex
structure compatible with $\Omega$ (see \cite{mei-sja}). Then we
denote $Q(M,\CL)$ simply by $Q(M)$. This is a canonical
finite-dimensional virtual representation $Q(M)$ of $K$, the
quantization of the symplectic manifold $M$. The spectrum of the
action of $\phi\in \k$ in $Q(M)$ should be the ``quantum" version
of the levels of energy of the Hamiltonian function $\langle
\mu,\phi\rangle $ on $M$ (see \cite{ver01} for survey). Guillemin
and Sternberg conjectured in 1982 that the multiplicity of the
irreducible representation $V_\xi$ of $K$ (of highest weight
$\xi\in \t^*_+\subset\k^*$) in the representation $Q(M)$ is equal
to $Q(M_\xi)$ and proved it for the case of K{\"a}hler manifolds. This
is summarized by the slogan: ``Quantization commutes with
Reduction". In other words, when $\xi=0$, we should have the
equality
$$\int_K \Tr_{Q(M)}(k)dk=\int_{M//K} {\rm ch}(\CL//K)
{\rm Todd}(M//K).$$

 Although
a fixed point formula exists for $\Tr_{Q(M)}(k)$, it is difficult
to extract the Guillemin-Sternberg conjecture directly from   the
Atiyah-Bott Lefschetz formula. Thus this conjecture (fundamental
for the credo of quantum mechanics) remained unproved for years.
Witten's inversion formula \cite{wit92}
$$\int_{\k}\left(\int_M e^{i\Omega(\phi)}p(\phi)\right)d\phi=
s_0(2i\pi)^{\dim \k} \vol(K)\int_{M//K}e^{i\Omega_{\rm red}}p_{\rm
red}$$ is in strong analogy with this conjecture. In particular,
apart from factors of $2i\pi$, the form $e^{i\Omega_{\rm red}}$ is
equal to ${\rm ch}(\CL//K)$. Meinrenken \cite{mei} used the
Atiyah-Bott Lefschetz formula  and symplectic cutting in a subtle
way  to give a proof of Guillemin-Sternberg conjecture for any
compact $K$-Hamiltonian manifold. This result was  extended
further  to singular symplectic quotients in Meinrenken-Sjamaar
\cite{mei-sja}.

\begin{definition}

Let $N$ be a prequantizable  Hamiltonian $K$-manifold with Kostant
line bundle $\CL$ such that the moment map is proper and the set
of zeroes of the Kirwan vector field $\kappa$ is compact. Define
$$Q(N):={\rm Index}(c_{\kappa,\CL}).$$
\end{definition}

Thus $Q(N)$ is a Fourier series of characters $\Tr(V_\xi)$.

\bigskip

{\bf Conjecture:}
   The
 multiplicity $m_\xi$ of the irreducible representation $V_\xi$
 in $Q(N)$  is equal to $Q(N_\xi)$.

When $N$ is compact, this is  the Guillemin-Sternberg conjecture.

\bigskip
Paradan \cite{par3}  proved this conjecture (in the spin context)
when $N:=G\lambda$ is an admissible regular elliptic coadjoint
orbit of a reductive real Lie group $G$ and $K$ the maximal
compact subgroup of $G$. Together with Theorem \ref{discrete},
this implies that irreducible representations $\Theta^K(\xi)$ (of
highest weight $\xi-\rho_{\k}$) of $K$ occurring in
Harish-Chandra's discrete series $\Theta^G(\lambda)|_K$ are such
that $\xi$ lies in the interior of the Kirwan polytope
$\mu(N)\cap\t^*_+$. This is a strong constraint on representations
appearing in $\Theta^G(\lambda)|_K$.

\begin{figure}[!h]
\begin{center}
\includegraphics{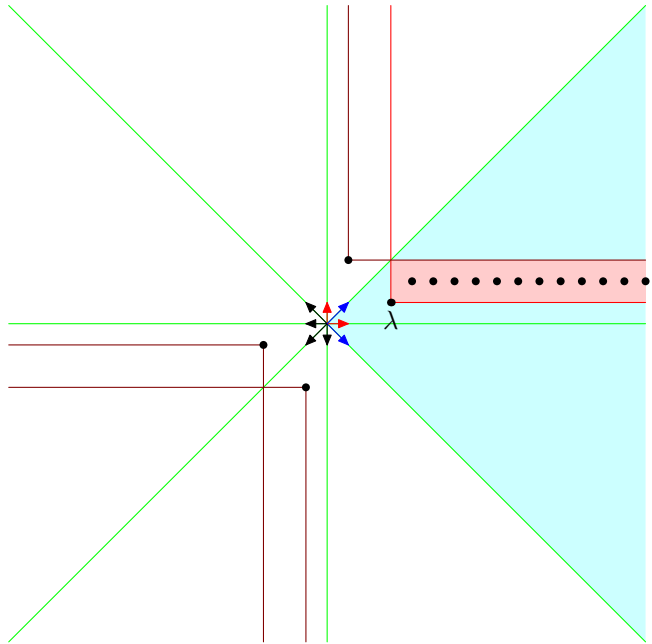}

\end{center}
\caption{ Restriction of Discrete Series and the Kirwan
polytope.}\label{kirS}
\end{figure}

\begin{example}\label{SO5}
Figure \ref{kirS} is the drawing for the restriction of the
representation $\Theta^G(\lambda)$ of $SO(4,1)$ to $SO(4)$. The
black dots are the $\xi$ such that $\Theta^K(\xi)$ occurs in
$\Theta^G(\lambda)$ (they all occur with multiplicity $1$). The
horizontal strip is the Kirwan polytope  $\mu(G\lambda)\cap
\t^*_+$.
\end{example}

\section{Arrangement of hyperplanes}\label{arrangement}

\subsection{Convolution of Heaviside distributions and cycles in the
complement of a set of hyperplanes}\label{Heaviside}

Let us consider a set $\CB:=\{\beta_1,\ldots,\beta_n\}$ of linear
forms $\beta_a$ on a vector space $V$ of dimension $r$, all in an
open half-space of $V^*$. We assume that the set $\CB$ spans
$V^*$. By definition, an element $\xi\in V^*$ is regular if it
does not lie in a cone spanned by $(r-1)$ elements of $\CB$. A
connected component of the set of regular elements is called a
chamber.

The convolution $H$ of the Heaviside distributions of  the
half-lines $\R^+\beta_a$ is a multivariate spline function on
$V^*$, that is, a locally polynomial function continuous on the
cone ${\rm Cone}(\CB)$ spanned by $\CB$.  Our problem is to
compute $H(\xi)$ at a particular point $\xi\in V^*$. In principle,
$H(\xi)$ is given by the following limit of integrals (on the
non-compact ``cycle" $V$ of dimension $r$, and in the sense of
Fourier transforms):

$$H(\xi)=\lim_{\epsilon\to 0}(2i\pi)^{-r}\int_{ V}e^{-i\langle \xi,v\rangle }
\frac{1}{\prod_{a=1}^n \langle \beta_a,v+i\epsilon\rangle }dv$$
where $\epsilon$ is in the dual cone to ${\rm Cone}(\CB)$.

Consider the complement of the hyperplanes defined by $\CB$ in the
complexified space $V_\C$:

$$V(\CB):=\{v\in V_\C; \langle
v,\beta\rangle \neq 0 \, \text{for all}\,  \beta\in \CB\}.$$

Jeffrey and Kirwan \cite{jef-kir95} introduced a residue calculus
on the space of functions defined on $V(\CB)$. A rational function
on $V(\CB)$ is of the form $R(v)=\frac{L(v)}{\prod_{a=1}^n \langle
\beta_a,v\rangle ^{n_a}}$ where $L(v)$ is a polynomial. The
following theorem  results from Jeffrey-Kirwan ideas, further
refined in \cite{Bri-Ver99} and \cite{sze-ver04}. We still denote
by $dv$ the holomorphic $r$-form $dv_1\wedge \cdots\wedge dv_r$ on
$V_\C$.
\begin{theorem}\label{cycle}\label{JK}
Let $\gc\subset {\rm Cone}(\CB)$ be a chamber. There exists a
compact oriented cycle $Z(\gc)$ of dimension $r$ contained in
$V(\CB)$ such that for any rational function $R$ on $V(\CB)$ and
any $\xi\in \gc$
$$\lim_{\epsilon\to 0}\int_{
V} e^{-i\langle \xi,v\rangle } R(v+i\epsilon)dv=\int_{Z(\gc)}
e^{-i\langle \xi,v\rangle } R(v)dv.$$
\end{theorem}

We gave in \cite{sze-ver04} a representative for the
$r$-dimensional cycle $Z(\gc)$ in $\C^r$ as the set of solutions
of $r$ real analytic equations related to  quantum cohomology.
Furthermore, we gave a simple  algorithm, further simplified   by
De Concini-Procesi \cite{con-pro}, to compute the homology class
of $Z(\gc)$ as a disjoint union of tori, so that integration on
$Z(\gc)$ is simply the algebraic operation of taking ordinary
iterated residues. Indeed, if $T({\beps})\subset V(\CB)$ is a
compact torus of the form in some coordinates
$(v_1,v_2,\ldots,v_r)\in \C^r:=V_\C$
$$T({\beps}):=\{v\in V(\CB); |v_k|=\epsilon_k,\text{ for }
k=1,\ldots, r\},$$ with $\beps:=[\epsilon_1\ll \epsilon_2\ll
\cdots \ll \epsilon_r]$  a sequence of increasing real numbers
(here $\epsilon_1\ll\epsilon_2$ meaning  that $\epsilon_2$ is
significantly greater than $\epsilon_1$, see \cite{sze-ver04} for
precise definitions), then the integration on $T({\beps})$ of a
function $F(v_1,v_2,\ldots, v_r)$ with poles on the hyperplanes
defined by  $\CB$ is
$$\frac{1}{(2i\pi)^r}\int_{T(\beps)} F(v_1,v_2,\ldots, v_r)dv=
\res_{v_r=0}\res_{v_{r-1}=0}\cdots \res_{v_1=0}
F(v_1,v_2,\ldots,v_r),$$ where each residue is taken assuming that
the variables with higher indices have a fixed, non-zero value.

Let me explain why this algorithm is efficient for computing the
convolution $H(\xi)$ of a large number of Heaviside distributions
in a vector space of small dimension. The usual way to compute
$H(\xi)$ is by induction on the cardinal of $\CB$. Here we fix
$\xi$ in a chamber $\gc$ and we  compute the cycle $Z(\gc)$
(depending on $\gc$) by induction on the dimension of $V$. It can
be done quite quickly using the maximal nested sets of De
Concini-Procesi, at least for classical root systems
\cite{bal-bec-coc-ver}.

\subsection{Intersection numbers on toric manifolds}\label{toric}
Let $T$ be a torus of dimension $r$  acting diagonally on
$N:=\C^n$ with weights $\CB:=[\beta_1,\beta_2,\ldots,\beta_n]$. We
assume that the cone ${\rm Cone}(\CB)$ spanned by the vectors
$\beta_a$ is an acute cone in $\t^*$ with non-empty interior. The
moment map $\mu:\C^n\to \t^*$ for the action of $T$ is
$\mu(z_1,\ldots,z_n)=\sum_{a=1}^n |z_a|^2\beta_a$. The reduced
space $N_\xi=\mu^{-1}(\xi)/T$ at a point $\xi\in {\rm Cone}(\CB)$
is {\em a toric variety}. It is an orbifold if $\xi$ is regular.
The space $N_\xi$ is still provided with a Hamiltonian action of
the full diagonal group $H:=(S^1)^n$ with Lie algebra
$\h:=\{\sum_{a=1}^n \nu_a J_a\}.$ The image of $N_\xi$ under the
moment map for $H$ is the convex polytope
$$P(\xi):=\{\sum_{a=1}^n x_a J^a\in \h^*; x_a\geq 0; \sum_{a=1}^n
x_a \beta_a=\xi\}.$$  Computing the volume of the polytope
$P(\xi)$ is the same as computing the symplectic volume of
$N_\xi$. All manifolds $N_\xi$ when $\xi$ varies in a chamber
$\gc$ are the same toric manifold $N_\gc$, the additional data
$\xi\in \gc$ being in one-to-one correspondence with the
Hamiltonian structure on $N_\gc$ coming from its identification
with the reduced space $N_\xi$.

The $T$-equivariant cohomology of $N$ is  $S(\t^*)$. For each
chamber $\gc$, the Kirwan map gives a surjective map
$\chi(p):=p_{\rm red}$ from $S(\t^*)$ to $\CH^*(N_\gc)$. The
following theorem
 allows us to compute integrals on toric manifolds.
\begin{theorem}\label{cycletoric}\cite{sze-ver04}
Let  $p\in S(\t^*)$, then
$$\int_{N_\gc}\chi(p)=(2i\pi)^{-r}\int_{Z(\gc)}\frac{p(\phi)}{\prod_{a=1}^n
\langle \beta_a,\phi\rangle}d\phi$$
\end{theorem}

Let $\xi\in \gc$ and let  $p(\phi):=\langle\phi,\xi\rangle$. Then
the cohomology class $p_{\rm red}$ is  the symplectic form of
$N_\gc$ determined by $\xi$. This way we obtain the formula:
\begin{corollary}\label{volumetoric}
Let $\xi\in \gc$, then
$$\vol(N_\xi)=\frac{1}{(2\pi)^r}\int_{Z(\gc)}
\frac{e^{-i\langle\xi,\phi\rangle}}{\prod_{a=1}^n
\langle\beta_a,\phi\rangle } d\phi.$$
\end{corollary}

We recall that the homology class of the cycle $Z(\gc)$ is
computed recursively  so that the preceding integral  is easily
calculated using iterated residues.

\section{Polytopes and Computations}\label{poly} It is well known
that many  theorems on toric varieties  have analogues in the
world of polytopes. With Brion, Szenes, Baldoni, Berline, we
carefully gave elementary proofs of the corresponding theorems on
polytopes, even if our inspiration came from equivariant
cohomology on Hamiltonian manifolds.

Let $\CB:=[\beta_1,\ldots,\beta_n]$ be a sequence of linear forms
on a vector space $V$ of dimension $r$ strictly contained in a
half-space of $V^*$. If $\xi\in V^*$, the partition polytope is
$$P_{\CB}(\xi):=\{\bx=[x_1,x_2,\ldots,x_n]\in \R^n; x_a\geq 0;
\sum_{a=1}^nx_a \beta_a=\xi\}.$$ Any polytope can be realized as a
partition polytope.

\begin{example} {\bf Transportation
polytopes.} Consider two sequences $[r_1,r_2,\ldots, r_k],$
$[c_1,c_2,\ldots, c_{\ell}]$ of positive numbers  with $\sum_i
r_i=\sum_j c_j$. Then ${\rm Transport}(k, \ell, r,c)$ is the
polytope consisting of all real matrices with $k$ rows and $n$
columns, with non-negative entries, and with sums of entries in
row $i$ equal to $r_i$ and in column $j$ equal to $c_j$. This is a
special case of a network polytope (see \cite{bal-del-ver},
\cite{bal-ver}).
\end{example}

The volume of $P_{\CB}(\xi)$ is equal to the value at $\xi$ of the
convolution of the Heaviside distributions supported on the
half-lines $\R^+\beta_a$. This becomes computationally hard if
there is a large number of convolutions. The volume of ${\rm
Transport}(k, \ell, r,c)$ necessitates
 the convolution of $k
\ell$ Heaviside distributions  in a space of dimension $k+\ell-1$.
For example, Beck-Pixton \cite{bec-pix} could compute, on parallel
computers, the volume of ${\rm Transport}(k, \ell, r,c)$ for
$k=10,\ell=10$, for special values $r_i=c_j=1$ (thus convoluting
$100$ linear forms in a $19$ dimensional space)  in 17 years of
computation time (scaled on 1 Ghz processor).

\begin{theorem}\label{volume}
Let $\gc$ be a chamber of ${\rm Cone}(\CB)$  and let  $\xi\in
\overline\gc$. Then
$$\vol(P_{\CB}(\xi))=(2i\pi)^{-r}
\frac{1}{(n-r)!}\int_{Z(\gc)}\frac{\langle \xi,v\rangle
^{n-r}}{\prod_{a=1}^n \langle \beta_a,v\rangle }dv.$$
\end{theorem}

Using De Concini-Procesi recursive determination of $Z(\gc)$, this
formula is expressed as a specific  sum of iterated residues.

Assume the $\beta_a$ span a lattice $\Lambda$ in $V^*$, and that
$\xi$ is in $\Lambda$. The discrete analogue of the volume of
$P_{\CB}(\xi)$ is the number $N_\CB(\xi)$ of integral points in
the rational polytope $P_{\CB}(\xi)$. A fundamental result of
Barvinok \cite{bar} asserts that $N_\CB(\xi)$ can be computed in
polynomial time, when $n$ is fixed.

The function $N_{\CB}(\xi)$ associates to the vector $\xi$  the
number of ways to represent the vector $\xi$ as a sum of a certain
number of vectors $\beta_a$. This  is  called the vector-partition
function of $\CB$. There is also a formula \cite{sze-ver03} for
$N_{\CB}(\xi)$ as an integral  on  the cycle $Z(\gc)$. This
integral formula has interesting theoretical applications, such as
information on the jumps of the partition function from chamber to
chamber. For example, the appearance of the five linear factors in
$g(a,b)$ (Formula (\ref{factor}) of Section \ref{inverse}) follows
from \cite{sze-ver03}. However, except for relatively good systems
$\CB$, this formula does not allow polynomial time computations. A
program for the counting of number of points in any rational
polytope  following Barvinok's algorithm is done in  Latte
\cite{latte}. For   systems not too far from unimodularity, our
programs based  on  integration on $Z(\gc)$, that is, on iterated
residues, are more efficient. It leads to the fastest computation
of number of integral points in network polytopes
\cite{bal-del-ver}, Kostant partition functions, weight
multiplicities $c_\mu^{\lambda}$ and  tensor product
multiplicities $c_{\lambda,\mu}^{\nu}$ of classical Lie algebras
(the bit size of the weights $\lambda,\mu,\nu$ can be very large
\cite{bal-bec-coc-ver}, \cite{coc}).

Finally, let me describe the local Euler-Maclaurin formula which
was conjectured by Barvinok-Pommersheim \cite{bar-pom}. It was
after observing the analogy of this conjecture with the
localization theorem (Theorem \ref{parter}) that I fully realized
the beauty of this conjecture. Nicole Berline and I  proved it by
using elementary means, based on the study of some valuations on
rational cones in an Euclidean space,

Let $P$ be a convex polytope in $\R^d$. For the sake of simplicity
we assume  that $P$ has {\em integral vertices}. Let $\CF$ be the
set of faces of $P$. For each face $F$ of $P$, the transverse cone
of $P$ along $F$  is a cone of dimension equal to the codimension
of $F$.

\begin{theorem}\label{eml}(Local Euler-Maclaurin formula.). For each face $F$, there exists a constant coefficients differential operator $D_F$
 (of infinite order), depending only on the transverse cone
 of $P$ along $F$, such that, for any polynomial function $\Phi$ on $\R^d$,
$$\sum_{\xi\in P\cap \Z^d}\Phi(\xi)=\sum_{F\in \CF} \int_FD_F
 (\Phi).$$
\end{theorem}

The detailed statement for any rational convex polytope and what
we really mean by ``depending only on" is in \cite{ber-ver06}.

The operators $D_F$ have rational coefficients and can be computed
in polynomial time when $d$  and the order of the expansion are
fixed, with the help of the Barvinok signed decomposition of cones
and the LLL short vector algorithm. The local property of $D_F$
means that if two polytopes $P$ and $P'$  are the same in a
neighborhood of a generic point of $F$, then the operators $D_F$
for $P$ and  $P'$ coincide.

The local Euler-Maclaurin formula gives in particular a local
formula for the number of integral points in $P$ or in the dilated
polytopes $tP$. The Ehrhart polynomial  $E(P)(t)$  is defined as
the number of integral points in $tP$, for $t$ a non-negative
integer. Then $E(P)(t)=\sum_{i=0}^n e_i t^{n-i}$, with
$e_0=\vol(P)$. Barvinok \cite{barsimplice} recently showed that
the (periodic) coefficients
 $e_i$ with $i\leq k$ can be computed in polynomial time, when $P$ is a rational simplex.
 We  hope to implement soon another  polynomial time algorithm for the
same problem based on our local formula.

Even though Time often prevails, in numerical computations as in
life, it was rewarding for us to see that our theoretical results
could help in effective computations.

\end{document}